\documentclass[12pt]{amsart}
\usepackage[dvips]{epsfig}
\usepackage{graphics}
\usepackage{latexsym}
\usepackage{verbatim}
\usepackage{amsmath}
\usepackage{amsthm}
\usepackage{amssymb}
\newtheorem{theorem}{Theorem}[section]

\newtheorem{proposition}[theorem]{Proposition}

\def\Remarks{\medskip\noindent{\bf Remarks: }}

\newcommand{\ens}[1]{\mathbb{#1}}

\newcommand{\N}{\mathbb{N}}

\newcommand{\R}{\mathbb{R}}

\def\derpar#1#2{\frac{\partial#1}{\partial#2}}

\begin{document}

\title[Large time behavior for soft potentials]
{Large time Behavior of the {\em a priori} bounds for the solutions to the spatially homogeneous 
Boltzmann equations with soft potentials}

\author{Laurent Desvillettes, Cl\'ement Mouhot}

\hyphenation{bounda-ry rea-so-na-ble be-ha-vior pro-per-ties
cha-rac-te-ris-tic}

\maketitle

\begin{abstract} 
We consider the spatially homogeneous Boltzmann equation for 
regularized soft potentials and Grad's angular cutoff. 
We prove that uniform (in time) bounds in 
$L^1((1+|v|^s)dv)$ and $H^k$ norms, $s,k \ge 0$ hold for its solution.
The proof is based on the mixture of estimates 
of polynomial growth in time of those norms 
together with the quantitative results of relaxation to equilibrium in $L^1$ obtained 
by the so-called ``entropy-entropy production'' method in the context of 
dissipative systems with slowly growing {\em a priori} bounds~\cite{ToscVill:cveq:2000}. 
\end{abstract}

\bigskip



\tableofcontents

\section{Introduction}
\setcounter{equation}{0}

This note is devoted to the study of the asymptotic behavior 
of solutions to the spatially homogeneous Boltzmann equation 
in the case of regularized soft potentials with Grad's angular cutoff. 
\smallskip

More precisely, we are concerned with the evolution of 
suitable norms which measure the asymptotic 
tail behavior (when $|v| \to +\infty$) of the distribution, and its smoothness. 
We shall prove bounds on the $L^1((1+|v|^q)dv)$ 
moments (resp. $H^k$ norms) of the 
distribution which are uniform with respect to time, 
provided that the initial datum belongs to 
$L^1((1+|v|^{q_0})dv) \cap H^{k_0}$ with $q_0,k_0$ big enough.
\medskip

The Boltzmann equation (Cf. \cite{Ce88} and \cite{CIP}) 
describes the behavior of a dilute gas when the only 
interactions taken into account are binary collisions.
In the case when the distribution function is
assumed to be independent on the position $x$, we obtain 
the so-called spatially homogeneous Boltzmann equation, which reads
 \begin{equation}\label{e1}
 \derpar{f}{t}(t,v)  = Q(f,f)(t,v), \qquad  v \in \R^N, \quad t \geq 0,
 \end{equation}
where $N \ge 2$ is the dimension. 
In equation~\eqref{e1}, $Q$ is the quadratic
Boltzmann collision operator, defined by the bilinear form
 \begin{equation*}\label{eq:collop}
 Q(g,f)(v) = \int _{\R^N \times \ens{S}^{N-1}} B(|v-v_*|, \cos \theta)
         \left(g'_* f' - g_* f\right) \, dv_* \, d\sigma,
 \end{equation*}
where we have used the shorthands $f=f(v)$, $f'=f(v')$, $g_*=g(v_*)$ and
$g'_*=g(v'_*)$. Moreover, $v'$ and $v'_*$ are parametrized by
 \begin{equation*}\label{eq:rel:vit}
 v' = \frac{v+v_*}2 + \frac{|v-v_*|}2 \, \sigma, \qquad 
 v'_* = \frac{v+v_*}2 - \frac{|v-v_*|}2 \, \sigma. 
 \end{equation*}
Finally, $\theta\in [0,\pi]$ is the deviation angle between 
$v'-v'_*$ and $v-v_*$ defined by $\cos \theta = (v'-v'_*)\cdot(v-v_*)/|v-v_*|^2$, 
and $B$ is the Boltzmann 
collision kernel determined by physics 
(related to the cross-section $\Sigma(v-v_*,\sigma)$ 
by the formula $B=|v-v_*| \, \Sigma$).  We also denote  
$$ Q^+(g,f)(v) =  \int _{\R^N \times \ens{S}^{N-1}} B(|v-v_*|, \cos \theta)\, g'_* f'\, dv_* \, d\sigma$$
the positive part of $Q$, and 
$$ L(g)(v) = \int _{\R^N \times \ens{S}^{N-1}} B(|v-v_*|, \cos \theta)\, g_*\, dv_* \, d\sigma $$
the linear operator appearing in the loss part of $Q$.
\medskip

Boltzmann's collision operator has the fundamental properties of
conserving mass, momentum and energy
  \begin{equation}
  \int_{\R^N}Q(f,f) \, \phi(v)\,dv = 0, \qquad
  \phi(v)=1,v,|v|^2, \label{CON}
  \end{equation}
and satisfying Boltzmann's $H$ theorem, which writes (at the formal level) 
  \begin{equation*} 
  - \frac{d}{dt} \int_{\R^N} f \log f \, dv = - \int_{\R^N} Q(f,f)\log(f) \, dv \geq 0.
  \end{equation*}
Boltzmann's $H$ theorem implies that (when  $B>0$ a.e.)
any equilibrium distribution function has the form of a Maxwellian distribution
  \begin{equation*}
  M(\rho,u,T)(v)=\frac{\rho}{(2\pi T)^{N/2}}
  \exp \left( - \frac{\vert u - v \vert^2} {2T} \right), 
  \end{equation*}
where $\rho\ge 0,\,u \in \R^N,\,T>0$ are the density, mean velocity
and temperature of the gas, defined by
  \begin{equation*}
  \rho = \int_{\R^N}f(v) \, dv, \quad u =
  \frac{1}{\rho}\int_{\R^N}vf(v) \, dv, \quad T = {1\over{N\rho}}
  \int_{\R^N}\vert u - v \vert^2f(v) \, dv, 
  \end{equation*}
and determined by the mass, momentum and energy of the initial datum thanks 
to the conservation properties (\ref{CON}). As a result of the process of entropy production 
pushing towards local equilibrium combined with the constraints (\ref{CON}), 
solutions are expected to converge to a unique Maxwellian equilibrium.
\medskip
 
This suggests for uniform bounds in time on the decay (in the $v$ variable)  and smoothness of the 
distribution $f=f(t,v)$. The main idea of this paper is to quantify this idea in a situation where
the uniform bounds are not obvious~: for so-called soft potentials.  
\medskip

More precisely, we shall consider the following assumptions on the collision kernel $B$: 
 \begin{itemize}
 \item[(H1)] It takes the following tensorial form (with $\Phi,b$ nonnegative functions)
	\[ B(|v-v_{*}|, \cos \theta) = \Phi(|v-v_{*}|) \, b(\cos \theta). \]

 \item[(H2)] The kinetic part $\Phi$ is $C^{\infty}$ and satisfies the bounds 
   \[ \forall \, z \in \R^N, \quad c_\Phi \, (1+|z|)^{\gamma} \le \Phi(|z|) 
      \le C_\Phi \, (1+|z|)^{\gamma}, \]
   \[ \forall \, z \in \R^N, \ p\in \N^*,  \quad |\Phi^{(p)} (|z|) | \le C_{\Phi,p}, \]
 with $\gamma \in (-2,0]$, and $c_\Phi, C_\Phi, C_{\Phi,p} >0$. 

\smallskip
 \item[(H3)] The angular part $\sigma \mapsto b(u \cdot \sigma)$ is integrable on 
 $\ens{S}^{N-1}$, and it satisfies the bound from below 
   \[ \forall \, \theta \in [0,\pi], \quad b(\cos \theta) \ge b_0 \]
 for some constant $b_0 >0$. 
 \end{itemize}
\smallskip

This includes the so-called ``mollified'' soft potentials with 
Grad's angular cutoff assumption (the word ``mollified'' is related to the 
singularity for small relative velocities). It does not include the very soft
potentials (that is the case when $\gamma \in (-N,-2]$). 
\smallskip

We shall systematically use the notations ($s\in\R$, $p\in [1, +\infty)$, $k \in \N$)
\[ \|f\|_{L^p_s}^p := \int_{\R^N} |f(v)|^p\, (1+ |v|^2)^{ps/2}\, dv, \]
and 
\[ \|f\|^2 _{H^k_s} := \sum_{0\le |i| \le k} \|\partial^i f\|_{L^2_s}^2, \]
where $\partial^i$ denotes the partial derivative related to the 
multi-index $i$.
\par
The Cauchy theory for equation~\eqref{e1} under assumptions (H1)-(H2)-(H3) 
is already known and is particularly simple (the collision operator is bounded). 
Using the arguments of Arkeryd~\cite{Ark72}, one can construct global nonnegative 
solutions in $L^1_2$. Uniqueness (in this class) follows from the boundedness of the operator 
(as a bilinear function in $L^1_2$). 
\medskip

As far as hard potentials (that is, $\gamma \in (0,1]$) or Maxwell molecules 
(that is, $\gamma =0$) are concerned, the propagation 
of the $L^1$ moments (that is, the $L^1_s$ norms for $s>2$) 
was proven in~\cite{elm} and~\cite{IkTr:74}. Moreover, the bounds 
were shown there to be uniform with respect to time. 
It was later noticed that for hard potentials 
those moments appear even if they don't initially exist,
under reasonable assumptions
(Cf. \cite{Desv:93}, and the improvements in~\cite{Wenn:momt:94,Wenn:momt:97,MiWe:99}). 
\par
Still for hard potentials (with angular cutoff), 
uniform in time estimates of $L^p$ norms or $H^k$ norms were first obtained
 in~\cite{Gust:L^p:86,Gust:L^p:88} and ~\cite{Wenn:rado:94}, and later simplified 
and systematically studied in~\cite{MV04}.
\par
In the case of (mollified) soft potentials (with angular cutoff), 
polynomially growing bounds on the $L^1$ moments
were first  obtained in~\cite{Desv:93} and later extended to 
the case of the Landau equation in~\cite[Part~I, Appendix~B]{Vill:these}  
and~\cite{ToscVill:cveq:2000}.
Polynomially growing bounds  
on the $L^p$ norms were also obtained in~\cite{ToscVill:cveq:2000}. 
\par
This paper is devoted to the obtention of uniform in time bounds on $L^1$ moments and $H^k$ norm in the 
setting of (mollified) soft potentials (with angular cutoff), where only polynomially growing bounds exist, as we just explained.
\medskip

We now state our main result. 
\begin{theorem}\label{theo:main}
Let $s >2$ and $k \ge 0$ be given, together with an initial datum $0 \le f_{in}\in L^1_2(\R^N)$. 
We consider the unique solution $f(t,v)\ge 0$ in $L^1_2$ to  equation~\eqref{e1} under 
assumptions (H1)-(H2)-(H3). Then
  \begin{itemize}
  \item[(i)] there exists $q_{0} >0$ (depending on $s$, but not on $k$) such that if
  $f_{in} \in L^1 _{2 s} \cap L^2 _{q_{0}}$, 
  the associated solution $f=f(t,\cdot)$ satisfies 
    \[ \sup_{t \geq 0} \| f(t,\cdot) \|_{L^1 _s} \le C(s) \]
  for some explicit bound $C(s) >0$; 
  \item[(ii)] there is $s_0 >0$ and $k' \ge k$ (both depend on $k$) such that if
  $f_{in} \in L^1_{s_0} \cap H^{k'}$, 
  the associated solution $f=f(t,\cdot)$ satisfies  
    \[ \sup_{t \geq 0} \| f(t,\cdot) \|_{H^k} \leq C(k) \]
  for some explicit bound $C(k) >0$.
  \end{itemize}
\end{theorem}
	
\Remarks 

1. In both points (i) and (ii) of this theorem, 
the assumptions on the initial datum are most probably not optimal, and 
are likely to be relaxed, up to technical refinements in the 
proofs (for example, in point (i), the weighted  $L^2$ space can be replaced by some weighted  $L^p$
space  for any $p>1$). 
We do not try here to look for such optimal 
assumptions, since we are more interested in showing how to obtain the 
uniform bounds. Note however that the sole assumption 
$f_{in} \in \cap_{s>0} L^1 _s$ is probably not 
sufficient to propagate uniformly the $L^1 _s$ norm for $s>2$, and 
we conjecture that it may be possible to construct some counter-examples 
in the same spirit as those constructed in~\cite{BoCe:99} 
in order to disprove Cercignani's conjecture for Maxwell molecules interactions.  
\smallskip

2. We then note that the assumptions on the collision kernel can also certainly be relaxed. We conjecture
that all derivatives on the kinetic part of the cross section are not really needed (probably one is enough),
and that the angular part need not really be bounded below. However, our proof depends strongly on the angular
cutoff, and one would need original extra arguments to treat the non cutoff case. It also does not work
for very soft potentials (see the remark at the end of the proof). 
\smallskip

3. We think that our proof could be adapted to the Landau  kernel with soft potential without too many changes.
However, too soft potentials like the Coulomb potential might not be reachable.
\smallskip

4. When $f_{in}$ belongs to ${\mathcal S}(\R^N)$ 
the Schwartz space of rapidly decaying $C^\infty$ function, then $f(t,\cdot) \in {\mathcal S}(\R^N)$ 
and the corresponding seminorms are bounded uniformly with respect to time. 
This is obtained thanks to Sobolev inequalities and standard interpolations between 
$L^1_s$ and $H^k$.
 In particular, uniform bounds 
of the form 
  \[ \forall \, t \ge 0, \ \forall \, v \in \R^N, \quad f(t,v) \le C \, (1+|v|)^{-q} \]
are available.
\smallskip 

5. A rough calculation shows that for point~(i) of this theorem, $q_0 = 26$ is sufficient 
in the case when $N=3$ and $\gamma = -1$.
\medskip

\section{Proof of slowly increasing bounds}
\setcounter{equation}{0}

In this section, we recall results on the slowly increasing 
polynomial bounds on the moments and $L^p$ norms of the solutions of equation~\eqref{e1} 
from~\cite{Desv:93,ToscVill:cveq:2000}, and we extend them to deal with the $H^k$ norms. 
\medskip

Estimates of linear growth in time on the moments were obtained 
in~\cite{Desv:93} in the case $\gamma >-1$, and sketched in~\cite{Vill:mou:98} 
and ~\cite{ToscVill:cveq:2000} for $\gamma >-2$.
We give here a precise statement together with a short proof.
\medskip
 
  \begin{proposition}\label{prop:linear}
  Let $s >2$. Then for any initial datum $f_{in} \in L^1 _s$, 
  the unique associated solution $f=f(t,\cdot )$ to equation~\eqref{e1} 
  under assumptions (H1)-(H2)-(H3) satisfies the bounds 
    \[ \forall \, t \geq 0, \quad  \| f(t,\cdot) \|_{L^1 _s} 
          \leq C_{0}(s) \, (1+t), \]
  for some explicit constant $C_{0}(s) >0$ depending only on the mass and $L^1 _s$ norm of $f_{in}$. 
  \end{proposition}


\begin{proof}[Proof of Proposition~\ref{prop:linear}]
We compute the time derivative of the $s$-th $L^1$ moment of $f$ 
thanks to the pre-postcollisional change of variable (see~\cite[Chapter~1, Section~4.5]{Vill:hand}):
 \[ \frac{d}{dt}  \| f(t,\cdot) \|_{L^1 _s} 
	    = \int_{\R^N \times \R^N \times \ens{S}^{N-1}} f \, f_* \, b \, 
	          \left( |v'_{*}|^s + |v'|^s - |v_{*}|^s - |v|^s \right) \, 
    \Phi(|v-v_{*}|) \, dv \, dv_{*} \, d\sigma. \]
Using then Povzner's inequality (Cf. \cite{Wenn:momt:97} for instance), we get (for some $C_+, K_->0$)
 \begin{multline*} 
  \int_{\ens{S}^{N-1}} \left( |v'_{*}|^s + |v'|^s - |v_{*}|^s - |v|^s \right) \, 
  b(\cos \theta) \, d\sigma \\  
  \le C_+ \, \left( |v|^{s-2} |v_{*}|^2 + |v_{*}|^{s-2} |v|^2\right)  
  - K_- \, \left(  |v|^s +  |v_{*}|^s \right). 
  \end{multline*}
Hence, using assumption (H2), for some $K_0>0$, 
\begin{equation}
\frac{d}{dt}  \| f(t,\cdot) \|_{L^1 _{s}} \le C_1 + C_2 \, \| f(t,\cdot) \|_{L^1 _{s-2}} 
                  - K_0 \, \| f(t,\cdot) \|_{L^1 _{s+\gamma}}. 
\end{equation}
We conclude by using an interpolation of $\| f(t,\cdot) \|_{L^1 _{s+ \gamma}}$ 
between $\| f(t,\cdot) \|_{L^1}$ and $\| f(t,\cdot) \|_{L^1 _{s-2}}$.
We see that
\begin{equation}
\frac{d}{dt} \| f(t,\cdot) \|_{L^1 _{s}} \le C_3(s),
\end{equation}
 so that 
\begin{equation}
\| f(t,\cdot) \|_{L^1 _{s}} \le C_4(s) \, (1 + t).
\end{equation}
\end{proof}

We now take care of the smoothness. The following result is a straightforward consequence 
of~\cite[Corollary~9.1]{ToscVill:cveq:2000} and general methods developed 
in~\cite{MV04}. It  essentially says  that the control of the regularity in our context 
can be obtained by the control of the moments. 

  \begin{proposition}\label{prop:reglent}
  Let $1<p<+\infty$ (resp. $k \in \N^*$). Let us consider $0 \le f_{in}\in L^1_2$ an initial datum 
  and $f=f(t,\cdot)$ the unique associated 
  solution to equation~\eqref{e1} under assumptions (H1)-(H2)-(H3). Then, 
  there are $C,s,\alpha >0$ depending on $p$ (resp. $C',s',\alpha' >0$ 
  depending on $k$) such that the following {\em a priori} estimates hold
    \begin{equation*}
    \left\{ 
    \begin{array}{l}\displaystyle
    \frac{d}{dt} \|f(t,\cdot)\|_{L^p} \le C \, \|f(t,\cdot)\|_{L^1 _s} ^\alpha, \vspace{0.3cm} \\ \displaystyle
    \frac{d}{dt} \|f(t,\cdot)\|_{H^k} \le C' \, \|f(t,\cdot)\|_{L^2 _{s'}} ^{\alpha'}. 
    \end{array}
    \right.
    \end{equation*}
  \end{proposition}

\begin{proof}[Proof of Proposition~\ref{prop:reglent}]
Concerning the first {\em a priori} bound on the $L^p$ norm, 
it is proven in~\cite[Proposition~9]{ToscVill:cveq:2000}. The proof 
is based on the regularity property of the gain part $Q^+$ of the collision 
operator in the following form (see~\cite{Wenn:rado:94,MV04}) 
  \[ \|Q^+(g,f)\|_{H^{(N-1)/2} _s} \le C \big( \|f\|_{L^1 _{1+2s}} \|g\|_{L^2 _{1+s}} 
                                            + \|f\|_{L^2 _{1+s}} \|g\|_{L^1 _{1+2s}} \big). \]
Then using that (for any derivative $\partial$)
  \[ \partial Q^+(g,f) = Q^+(\partial g, f) + Q^+(g, \partial f) \]
thanks to the translation invariance, to Cauchy-Schwartz type inequalities like 
  \[ \|f\|_{L^1 _s} \le C \, \|f\|_{L^2 _{s+q}} \]
for some $C,q>0$, and to some classical interpolation in the $H^k$ spaces, we deduce that 
  \begin{equation} \label{eq:regQ^+}
  \|Q^+(f,f)\|_{H^{k+(N-1)/2} _s} \le C \|f\|_{L^2 _{s+w}} \|f\|_{H^k _{s+w}}
  \end{equation}
for any $s,k \ge 0$ and some $C,w>0$. Now let us consider the time 
derivative of the square of the $L^2$ norm of $\partial^k f$ for some multi-index $k$ 
with $|k| \ge (N-1)/2$. We get 
  \begin{multline*}
  \frac{d}{dt} \|\partial^k f\|_{L^2} ^2 \le 
  C\, \|Q^+(f,f)\|_{H^k _{-\gamma}} \|\partial^k f\|_{L^2 _\gamma} \\ + C \, 
  \left( \sum_{0<l\le k} \big\| (\partial^l L(f)) (\partial^{k-l} f) \big\|_{L^2 _{-\gamma}} \right) 
   \|\partial^k f\|_{L^2 _\gamma} 
  - K \, \|\partial^k f\|^2 _{L^2 _{\gamma}}, 
  \end{multline*}
where we have $L(f)= C_b \, (\Phi * f)$ ($C_b$ is the $L^1$ norm of 
$b$ on the sphere $\ens{S}^{N-1}$). The back term comes from 
the classical lower bound 
  \[ L(f) \ge K \, (1+|v|)^{\gamma} \]
for some constant $K>0$ depending on the mass and entropy of the initial datum 
(see~\cite{Ark72} for instance). 

Then on the one hand equation~\eqref{eq:regQ^+} yields 
  \[ \|Q^+(f,f)\|_{H^k _{-\gamma}} \le C_+ \, \|f\|_{L^2 _w} \|f\|_{H^{k-(N-1)/2} _w} \]
for some explicit constants $C_+,w>0$, and then  by interpolation
  \[ \|Q^+(f,f)\|_{H^k _{-\gamma}} \le C_+ \, \|f\|_{L^2 _{w_+}} ^{1+2\theta_+} 
      \|f\|_{H^k _\gamma}^{1-2\theta_+} \]
for some explicit constants $C_+,w_+,\theta_+ >0$. On the other hand the convolution 
structure of $L(f)$ together with the smoothness assumption on $\Phi$ in (H2) yields 
easily 
  \[ \left( \sum_{0<l\le k} \big\| (\partial^l L(f)) (\partial^{k-l} f) \big\|_{L^2 _{-\gamma}} \right) 
      \le C \, \| f \|_{L^2_w} \|f\|_{H^{k-1}_w} \]
for some constants $C,w>0$ and thus by interpolation 
  \[ \left( \sum_{0<|l|\le k} \big\| (\partial^l L(f)) (\partial^{k-l} f) \big\|_{L^2 _{-\gamma}} \right) 
      \le C_L \, \|f\|_{L^2 _{w_L}} ^{1+2\theta_L} 
      \|f\|_{H^k _\gamma}^{1-2\theta_L} \]
for some explicit constants $C_L,w_L,\theta_L >0$ (depending on $k$). 
Thus if $\theta_0 = \min\{\theta_+,\theta_L\}$,  
we easily obtain for some $\bar w >0$, 
  \[ \frac{d}{dt} \| f\|_{H^k}^2 \le C_0 \, \|f\|_{L^2 _{\bar w}} ^{1+2\theta_0} 
     \|f\|_{H^k_{\gamma}}^{2-2\theta_0}  - K \, \|f\|^2_{H^k_{\gamma}}. \]
Finally, we use the inequality 
  \[ \forall \, X \ge 0, \quad A \, X^{1-\delta} - K X \le C \, A^{1/\delta} \] 
(with
$C\equiv C(K, \delta)$) to conclude the proof.
\end{proof}
%

\section{Proof of uniform bounds}
\setcounter{equation}{0}

In this section, we combine the results of Section 2 
with the quantitative results of convergence to equilibrium 
obtained in~\cite{ToscVill:cveq:2000}. We 
conclude in this way the proof of Theorem~\ref{theo:main}. 
\medskip

Let us recall the quantitative result of trend to equilibrium we shall use. 
We denote by $M= M(\rho, u, T)$ the Maxwellian with parameters 
$\rho, u,T$ corresponding to the initial datum.

  \begin{proposition}\label{prop:eep}
  Let us consider an initial datum $0 \le f_{in} \in L^1 _2$ and $\tau>0$.   
  Then there exists $q_{0}>0$
  such that if $f_{in} \in L^2 _{q_{0}}$, 
  the unique associated solution $f=f(t,\cdot )$ of equation~\eqref{e1} under assumptions 
  (H1)-(H2)-(H3) satisfies 
    \[ \forall \, t \geq 0, \quad  \| f(t,\cdot) - M \|_{L^1} \leq C_{1} \, (1+t)^{-\tau}  \]
  for some explicit bound $C_{1} >0$ depending only on $\tau$, $\rho$, and 
  the $L^2_{q_{0}}$ norm of $f_{in}$.
  \end{proposition} 

\begin{proof}[Proof of Proposition~\ref{prop:eep}]
This result is a particular case of more general results 
in~\cite{ToscVill:cveq:2000} (see Proposition~6 in this paper). 
Indeed~\cite[Theorem~11]{ToscVill:cveq:2000} implies the conclusion 
of Proposition~\ref{prop:eep} as soon as $f_{in}$ satisfies a lower bound 
of the form $f_{in} \ge K_0 \, e^{-A_0 |v|^2}$. This assumption 
can be relaxed thanks to~\cite[Theorem~5.1]{Mbinf}, which shows that 
this lower bound appears immediately under the assumption we have 
on the initial datum (in particular the assumption of finite 
entropy for~\cite[Theorem~5.1]{Mbinf} is implied by $f_{in} \in L^2 _{q_{0}}$). 
\end{proof}


Now we can conclude the proof of Theorem~\ref{theo:main} by gathering 
this proposition with the results of Section 2. 
\bigskip

{\bf Proof of point~(i)}: Assume that $f_{in} \in L^1 _{2s} \cap L^2 _{q_{0}}$. 
On the one hand, from Proposition~\ref{prop:linear}, the unique associated solution satisfies  
  \[ \forall \, t \geq 0, \quad \| f(t,\cdot) \|_{L^1 _{2s}} \leq C_{0} \, (1+t). \]
On the other hand, from Proposition~\ref{prop:eep} (with $\tau=1$), it satisfies 
  \[ \forall \, t \geq 0, \quad  \| f(t,\cdot) - M \|_{L^1} \leq C_{1} \, (1+t)^{-1}.  \]
We deduce that for any $t \geq 0$, 
  \begin{eqnarray*}
    \| f(t,\cdot) \|_{L^1 _{s}} &\le& \| M \|_{L^1 _{s}} 
    + \| f(t,\cdot) - M \|_{L^1 _{s}} \\  
    &\le& \| M \|_{L^1 _{s}} + \| f(t,\cdot) 
    - M \|_{L^1} ^{1/2} \| f(t,\cdot) - M \|_{L^1 _{2s}} ^{1/2} \\
    &\le& \| M \|_{L^1 _{s}} +  C_{1}  ^{1/2} (1+t)^{-1/2} \big( \| f(t,\cdot) \|_{L^1 _{2s}} 
    + \|M\|_{L^1 _{2s}} \big)^{1/2} \\
    &\le& \| M \|_{L^1 _{s}} +  C_{1}  ^{1/2} (1+t)^{-1/2} \big( C_{0} \, (1+t) 
    + \|M\|_{L^1 _{2s}} \big)^{1/2} \\
    &\le& C(s) < + \infty .
  \end{eqnarray*}
This concludes the proof of point (i).
\bigskip

{\bf Proof of point~(ii)}: First let us prove the uniform bound in the case $k=0$. 
In fact we shall prove uniform bounds on any $L^p$ norms, $1<p<+\infty$. 
From Proposition~\ref{prop:reglent} we have for any $p\in (1,+\infty)$
  \begin{equation}\label{yuyu}
 \frac{d}{dt} \|f(t,\cdot)\|_{L^p} \le C \, \|f(t,\cdot)\|_{L^1 _s} ^\alpha 
\end{equation}
for some explicit $C,s,\alpha >0$ (depending on $p$). 
\par
We assume enough $L^1$ moments bounded on the
initial datum, and enough derivatives in $L^2$. Then, thanks to Sobolev inequalities, the initial datum
is in $L^p$ with $p>2$. By standard interpolations, the initial datum has enough moments bounded in $L^2$. 
As a consequence, we can use point~(i), and obtain that 
$\|f(t,\cdot)\|_{L^1 _s}$ is uniformly bounded for 
all $t$. Using once again enough derivatives in $L^2$ of the initial datum and Sobolev inequalities, 
we get $L^p$ bounds  (for any $p\in ]1,+\infty[$) on the initial datum. Consequently, (\ref{yuyu})
yields 
  \begin{equation}\label{yu}
 \|f(t,\cdot)\|_{L^p} \le C_0(p) \, (1+t) 
\end{equation}
for some explicit constant $C_0(p)>0$.
\par
Then for any $p\in (1,+\infty)$ (using Proposition~3.1 and (\ref{yu}) for $2p$ instead of $p$) 
  \begin{eqnarray*}
    \| f(t,\cdot) \|_{L^p} &\le& \| M \|_{L^p} 
    + \| f(t,\cdot) - M \|_{L^p} \\  
    &\le& \| M \|_{L^p} + \| f(t,\cdot) 
    - M \|_{L^1} ^{1/(2p-1)} \| f(t,\cdot) - M \|_{L^{2p}} ^{1-1/(2p-1)} \\
    &\le& \| M \|_{L^p} +  C_{1}  ^{1/(2p-1)} (1+t)^{- 2(p-1)/(2p-1)} \big( \| f(t,\cdot) \|_{L^{2p}} 
    + \|M\|_{L^{2p}} \big)^{1-1/(2p-1)} \\
    &\le& \| M \|_{L^p} +  C_{1}  ^{1/(2p-1)} (1+t)^{- 2(p-1)/(2p-1)} \big( C_0(2p) \, (1+t) 
    + \|M\|_{L^{2p}} \big)^{1-1/(2p-1)} \\
    &\le& C(p) < + \infty.
  \end{eqnarray*}

Let us now assume that $k \ge 1$. From Proposition~\ref{prop:reglent}, we have 
  \[ \frac{d}{dt} \|f(t,\cdot)\|_{H^k} \le C' \, \|f(t,\cdot)\|_{L^2 _{s'}} ^{\alpha'} \]
for some explicit constants $C',s',\alpha'>0$ (depending on $k$). From the previous study, 
by assuming enough $L^1$ moments and $H^k$ bounds on the initial datum, we can assume that 
$\| f(t,\cdot) \|_{L^p}$ is uniformly bounded for all $t$. Using then point~(i) and a standard 
interpolation, we see that
$\|f(t,\cdot)\|_{L^2 _{s'}}$ is uniformly bounded for all $t$. 
Hence for any $k \ge 1$, we have 
  \begin{equation}\label{yo}
 \|f(t,\cdot)\|_{H^k} \le C_0(k) \, (1+t) 
\end{equation}
for some explicit constant $C_0(k)>0$.
\par
Then for any $k$, using Proposition~3.1 with $\tau=1$ and (\ref{yo}) with $2k + (N+1)/2$ instead of $k$ 
and the continuous embedding $L^1(\R^N) \hookrightarrow H^{-(N+1)/2}(\R^N)$, we have 
  \begin{eqnarray*}
    \| f(t,\cdot) \|_{H^k} &\le& \| M \|_{H^k} 
    + \| f(t,\cdot) - M \|_{H^k} \\  
    &\le& \| M \|_{H^k} + \| f(t,\cdot) 
    - M \|_{H^{-(N+1)/2}} ^{1/2} \| f(t,\cdot) - M \|_{H^{2k+(N+1)/2}} ^{1/2} \\
    &\le& \| M \|_{H^k} +  C_{1}  ^{1/2} (1+t)^{-1/2} \big( \| f(t,\cdot) \|_{H^{2k+(N+1)/2}} 
    + \|M\|_{H^{2k+(N+1)/2}} \big)^{1/2} \\
    &\le& \| M \|_{H^k} +  C_{1}  ^{1/2} (1+t)^{-1/2} \big( C_{0} \, (1+t) 
    + \|M\|_{H^{2k+(N+1)/2}} \big)^{1/2} \\
    &\le& C(k) < + \infty.
  \end{eqnarray*}
This concludes the proof. 
\bigskip

\Remarks  

1. Our analysis does not work for (mollified) very soft potentials. What happens then is that 
(if we denote by $m_s$ the $s$-th moment in $L^1$
of $f$), 
\[ \frac{d}{dt} m_s \le C_0 + C_1 \, m_{s-a} m_{a} - K \, m_{s+\gamma} \] 
for all $a \in [0,s]$, so that (in dimension $N=3$ with $-3<\gamma<-2$)
\[ m_s (t) \le C \, \left(1+t^{s/2 -1} \right) . \]
However, this estimate doesn't seem sufficient to obtain any rate of convergence to equilibrium.
A rough calculation shows that an estimate 
in $t^{\lambda s}$ instead of $t^{s/2}$ (with $\lambda<1/2$) 
could be the minimum required in order to get some rate of convergence to equilibrium 
with the ``entropy-entropy production'' method. 
Note however that for the Landau kernel for (mollified) very soft potentials  
(although not for the limiting Coulomb case) such estimates are available (see~\cite{ToscVill:cveq:2000}), 
suggesting that our method applies as well for this model.  
\smallskip

\smallskip

2. We conclude with a last remark: once bounds which are uniform in 
time have been  proven, they can be used in order to prove 
directly the rate of convergence toward equilibrium like in \cite{ToscVill:cveq:99} 
(that is, without entering the details of the method of ``slowly growing {\it a priori} estimates''
devised by G. Toscani and C. Villani in \cite{ToscVill:cveq:2000}). Note however that 
in order to get the bounds on moments which are uniform in time, this method 
(of ``slowly growing {\it a priori} estimates'') is used,
so it really seems unavoidable.

\bigskip
\noindent
{\bf{Acknowledgment}}: Support by the European network HYKE, funded by the EC as
contract HPRN-CT-2002-00282, is acknowledged.

\end{document}